\documentclass[12pt]{article}
\usepackage{epsfig}
\usepackage{amssymb}
\usepackage{graphicx}
\usepackage{amsmath}
\usepackage{amsfonts}
\usepackage{amsthm}
\usepackage{amscd}
\usepackage{arydshln}
\usepackage{mathtools}
\usepackage{units}
\usepackage{hyperref}
\usepackage[all,cmtip]{xy}
\usepackage{diagxy}
\usepackage{xymatrix}
\usepackage[all]{xy}
\usepackage{fancyhdr}
\usepackage[nottoc,numbib]{tocbibind}
\usepackage{graphicx}
\usepackage{mathrsfs}
\usepackage[perpage]{footmisc}
\usepackage{titlesec}
\usepackage{color}
\usepackage{geometry}
\usepackage{multicol}
\usepackage{slashed}
\usepackage{appendix}
\usepackage{changepage}

\newtheoremstyle{basic}{11pt}{11pt}{}{}{\bfseries}{.}{0.5em}{}
\newtheoremstyle{proof}{11pt}{11pt}{}{}{\scshape}{:}{0.5em}{}

\newtheorem{prop}{Proposition}
\newtheorem{cor}{Corollary}

\newtheorem{thm}{Theorem}
\newtheorem{que}{Question}

\theoremstyle{basic}

\theoremstyle{basic}

\theoremstyle{basic}

\theoremstyle{proof}
\newtheorem*{prf}{Proof}

\begin{document}

\pagestyle{plain}

\title{The Lie Algebra of S-unitary Matrices, Twisted Brackets and Quantum Channels}

\author{Clarisson Rizzie Canlubo}

\maketitle

\begin{abstract}
A dimension formula was given in \cite{caalim} in order to partially classify the Lie algebras of $S$-unitary type. The natural question of when $\mathfrak{u}_{S}$ and $\mathfrak{u}_{T}$ are isomorphic is left unanswered. In this article, we will give an answer to this question using the notion of quantum channels and their Kraus representation. In line with this, we will also discuss linearly twisted versions of the usual commutator bracket and its relation to the standard Lie algebra structure on $M_{n}(\mathbb{C})$. Finally, we will mention some problems that are still unanswered in relation to $S$-unitary type matrices and twisted brackets.
\end{abstract}

\section{Introduction} \label{intro}

Let $S\in M_{n}(\mathbb{C})$. Then, the subspace $\mathfrak{u}_{S}=\left\{ X\in M_{n}(\mathbb{C}) | SX^{\ast}=-XS \right\}$ is a Lie algebra with respect to the usual bracket of matrices, given as $[A,B]=AB-BA$ for any $A,B\in \mathfrak{u}_{S}$. The subset $U_{S}=\left\{ X\in M_{n}(\mathbb{C}) | SX^{\ast}=X^{-1}S \right\}$ of $M_{n}(\mathbb{C})$ is a Lie subgroup of $GL_{n}(\mathbb{C})$ whose Lie algebra is $\mathfrak{u}_{S}$.

In \cite{caalim}, the dimension of $\mathfrak{u}_{S}$ is given in terms of the spectral properties of $S$. Hence, if $S$ and $T$ are unitarily similar the Lie algebras $\mathfrak{u}_{S}$ and $\mathfrak{u}_{T}$ have the same dimension. Although this is not enough to conclude whether $\mathfrak{u}_{S}$ are isomorphic to $\mathfrak{u}_{T}$ as Lie algebras, this turns out to be the case according to the following proposition.

\begin{prop}\label{P1}
If $S$ and $T$ are unitarily similar then $\mathfrak{u}_{S}\cong\mathfrak{u}_{T}$ as Lie algebras.
\end{prop}

\begin{prf}
Suppose $S=V^{\ast}TV$ for some unitary $V$. Then, for any $X\in\mathfrak{u}_{S}$ we have

\[ V^{\ast}TVX^{\ast}=SX^{\ast}=-XS=-XV^{\ast}TV \]

\noindent and so, we have

\[ T(VXV^{\ast})^{\ast}=-(VXV^{\ast})T. \]

\noindent Thus, $\mathfrak{u}_{S}\stackrel{\phi_{V}}{\longrightarrow}\mathfrak{u}_{T}, X\mapsto VXV^{-1}$ gives the desired isomorphism. $\blacksquare$

\end{prf}

The isomorphism $\phi_{V}$ given in the proof of Proposition (\ref{P1}) turns out to be the most general one as indicated in the following theorem.

\begin{thm}\label{T1}
If $\mathfrak{u}_{S}\stackrel{\phi}{\longrightarrow}\mathfrak{u}_{T}$ is a Lie algebra isomorphism then $\phi(X)=VXV^{-1}$ for some invertible $V$.
\end{thm}

Also, a partial converse of Proposition (\ref{P1}) is a corollary of Theorem (\ref{P1}) as indicated in the next theorem.

\begin{thm}\label{T2}
If the Lie algebras $\mathfrak{u}_{S}$ and $\mathfrak{u}_{T}$ are isomorphic then the stabilizers of $S$ and $T$ under the conjugation action of $GL_{n}(\mathbb{C})$ on $M_{n}(\mathbb{C})$ are conjugate subgroups.
\end{thm}

\noindent  We will prove Theorems (\ref{T1}) and (\ref{T2}) in section [\ref{proof}]. Note that although the entries of the matrices in $\mathfrak{u}_{S}$ are complex numbers, the Lie algebra $\mathfrak{u}_{S}$ is strictly a \textit{real} Lie algebra. Whether $\mathfrak{u}_{S}$ is a complex Lie algebra depends on the existence of a complex structure $J$ (an endomorphism $J$ such that $J^{2}=-I$) which bilinearly commutes with $[,]$, i.e. $[J(X),Y]=J[X,Y]=[X,J(Y)]$ for any $X,Y\in \mathfrak{u}_{S}$.

\section{Twisted Lie Brackets}\label{twisted}

Using a linear map $M_{n}(\mathbb{C})\stackrel{\psi}{\longrightarrow}M_{n}(\mathbb{C})$, one can define a bilinear form $[,]_{\psi}$ on $M_{n}(\mathbb{C})$ as follows. For any $X,Y\in M_{n}(\mathbb{C})$, define $[,]_{\psi}$ by $[X,Y]_{\psi}=X\psi(Y)-Y\psi(X)$. Clearly, $[,]_{\psi}$ is skew-symmetric for any linear map $\psi$. In the event that $[,]_{\psi}$ defines a Lie bracket on $M_{n}(\mathbb{C})$, we will call ${,}_{\psi}$ the $\psi$-twisted Lie bracket on $M_{n}(\mathbb{C})$. However, the Jacobi identity is satisfied only for certain linear maps $\psi$. In general, any linear map $M_{n}(\mathbb{C})\stackrel{\psi}{\longrightarrow}M_{n}(\mathbb{C})$ takes the form $\psi(X)=\sum\limits_{i=1}^{m}A_{i}XB_{i}^{\ast}$ for some matrices $A_{i}, B_{i}\in M_{n}(\mathbb{C})$. If $B_{i}=I$ for all $i=1,\dots, m$ then $\psi(X)=AX$ for all $X\in M_{n}(\mathbb{C})$, where $A=\sum\limits_{i=1}^{m}A_{i}$. In this case, we have

\begin{eqnarray*}
\sum\limits_{\circlearrowleft} [X,[Y,Z]_{\psi}]_{\psi} &=& X\psi(Y\psi(Z))-X\psi(Z\psi(Y))-Y\psi(Z)\psi(X)+Z\psi(Y)\psi(X)\\
&+& Y\psi(Z\psi(X))-Y\psi(X\psi(Z))-Z\psi(X)\psi(Y)+X\psi(Z)\psi(Y)\\
&+& Z\psi(X\psi(Y))-Z\psi(Y\psi(X))-X\psi(Y)\psi(Z)+Y\psi(X)\psi(Z)\\
& & \\
&=& XAYAZ-XAZAY-YAZAX+ZAYAX\\
&+& YAZAX-YAXAZ-ZAXAY+XAZAY\\
&+& ZAXAY-ZAYAX-XAYAZ+YAXAZ\\
&=& 0
\end{eqnarray*}

\noindent where the leftmost sum indicates the cyclic sum over $X,Y$ and $Z$. Thus, we have proven the following proposition.

\begin{prop}
If $\psi(X)=AX$ for some $A\in M_{n}(\mathbb{C})$ then $[,]_{\psi}$ defines a Lie bracket on $M_{n}(\mathbb{C})$.
\end{prop}

The essential ideas of the Kraus representation of a linear map $M_{n}(\mathbb{C})\stackrel{\psi}{\longrightarrow}M_{n}(\mathbb{C})$ says that for any $X\in M_{n}(\mathbb{C})$, we have

\[ \psi(X)=\sum\limits_{i=1}^{m}A_{i}XB_{i}, \]

\noindent the case when the Choi matrix $J(\psi)$ of $\psi$ is of rank $m$. In the event that the rank of $J(\psi)$ is one, i.e. the sum above consists of only one summand, the following proposition gives a necessary condition when $[,]_{\psi}$ defines a Lie bracket on $M_{n}(\mathbb{C})$.

\begin{prop}
Let $\psi(X)=AXB$ for some $A,B\in M_{n}(\mathbb{C})$. Suppose the image of the map $M_{n}(\mathbb{C})\longrightarrow M_{n}(\mathbb{C}), X\mapsto AX$ is a subspace of $C(B)$, the commuting ring of $B$. Then $[,]_{\psi}$ defines a Lie bracket on $M_{n}(\mathbb{C})$. In particular, if $A$ is invertible then $\psi$ is of the form $\psi(X)=MX$ for some invertible matrix $M$.
\end{prop}

\begin{prf}
For any $X\in M_{n}(\mathbb{C})$, the assumption that $AX$ belongs to $C(B)$ implies that $(AX)B=B(AX)$. Thus, we have

\begin{eqnarray*}
\sum\limits_{\circlearrowleft} [X,[Y,Z]_{\psi}]_{\psi} &=& X\psi(Y\psi(Z))-X\psi(Z\psi(Y))-Y\psi(Z)\psi(X)+Z\psi(Y)\psi(X)\\
&+& Y\psi(Z\psi(X))-Y\psi(X\psi(Z))-Z\psi(X)\psi(Y)+X\psi(Z)\psi(Y)\\
&+& Z\psi(X\psi(Y))-Z\psi(Y\psi(X))-X\psi(Y)\psi(Z)+Y\psi(X)\psi(Z)\\
& & \\
&=& XAY(AZB)B-XAZ(AYB)B-Y(AZB)(AXB)\\
&+& Z(AYB)(AXB)+YAZ(AXB)B-YAX(AZB)B\\
&-& Z(AXB)(AYB)+X(AZB)(AYB)+ZAX(AYB)B\\
&-& ZAY(AXB)B-X(AYB)(AZB)+Y(AXB)(AZB)\\
& & \\
&=& XAY(AZB-BAZ)B-XAZ(AYB-BAY)B)\\
&+& YAZ(AXB-BAX)B-YAX(AZB-BAZ)B\\
&+& ZAX(AYB-BAY)B-ZAY(AXB-BAX)B\\
&=& 0
\end{eqnarray*}

\noindent for any $X,Y,Z\in M_{n}(\mathbb{C})$. This proves the first claim. Now, if $A$ is invertible then the image of $X\mapsto AX$ is $M_{n}(\mathbb{C})$. This implies that $B$ is necessarily a scalar matrix. Taking $M=BA$ proves the second claim. $\blacksquare$
\end{prf}

The canonical bracket on $M_{n}(\mathbb{C})$ restricts to a Lie bracket on the subspace $\mathfrak{u}_{S}$ described in section [\ref{intro}]. The following proposition gives a necessary condition when this is true for the twisted Lie bracket $[,]_{\psi}$.

\begin{prop}\label{P2}
Let $\psi(X)=AX$ for some $A\in M_{n}(\mathbb{C})$. Then, $[,]_{\psi}$ restricts to a Lie bracket on $\mathfrak{u}_{S}$ if and only if $A$ is $S$-Hermitian.
\end{prop}

\begin{prf}
Suppose $A$ is $S$-Hermitian. Then, for any $X,Y\in\mathfrak{u}_{S}$, we have

\begin{eqnarray*}
S[X,Y]_{\psi}^{\ast} &=& S(X\psi(Y)-Y\psi(X))^{\ast} \\
&=& SY^{\ast}A^{\ast}X^{\ast}-SX^{\ast}A^{\ast}Y^{\ast}\\
&=& -YSA^{\ast}X^{\ast}+XSA^{\ast}Y^{\ast}\\
&=& -YASX^{\ast}+XASY^{\ast}\\
&=& YAXS-XAYS\\
&=& (YAX-XAY)S\\
&=& -[X,Y]_{\psi}S
\end{eqnarray*}

\noindent Thus, $-[X,Y]_{\psi}\in \mathfrak{u}_{S}$, and so $[,]_{\psi}$ restricts to a Lie bracket on $\mathfrak{u}_{S}$.

Conversely, suppose $[,]_{\psi}$ restricts to a Lie bracket on $\mathfrak{u}_{S}$. Then, for any $X,Y\in \mathfrak{u}_{S}$, we have

\[ -YSA^{\ast}X^{\ast} + XSA^{\ast}Y^{\ast} = -YASX^{\ast} + XASY^{\ast} \]

\noindent from the above computation. Thus, we have

\[ Y(AS-SA^{\ast})X^{\ast} = X(AS-SA^{\ast})Y^{\ast} \]

\noindent for any $X,Y\in \mathfrak{u}_{S}$. Taking $Y=iI\in\mathfrak{u}_{S}$, we see that

\[ (AS-SA^{\ast})X^{\ast}=-X(AS-SA^{\ast}) \]

\noindent and taking $X=I$, we get $SA^{\ast}=AS$. $\blacksquare$

\end{prf}

It is a curiosity to know which linear maps $\psi$ whose brackets $[,]_{\psi}$ induce Lie algebra structures on $M_{n}(\mathbb{C})$ isomorphic to the canonical one. We partially answer this in the next proposition.

\begin{prop}\label{P5}
If the Lie bracket $[,]_{\psi}$, with $\psi(X)=AXB^{\ast}$ for all $X\in M_{n}(\mathbb{C})$, coincides with the canonical Lie bracket on $M_{n}(\mathbb{C})$ then $B^{\ast}=A^{-1}$.
\end{prop}

\begin{prf}
For any $X\in M_{n}(\mathbb{C})$ we have

\[ 0= [X,I] = [X,I]_{\psi} = X\psi(I)-I\psi(X). \]

\noindent And so, we have $\psi(X)=X\psi(I)$ for all $X\in M_{n}(\mathbb{C})$. Since $[I,X]=0$, we also have $\psi(X)=\psi(I)X$ for all $X\in M_{n}(\mathbb{C})$. Thus, $\psi(I)=AB^{\ast}$ is central. Thus, for any $X,Y\in M_{n}(\mathbb{C})$ we have

\[ [X,Y] = [X,Y]_{\psi} = X\psi(Y)-Y\psi(X) = X\psi(I)Y-Y\psi(I)X = \psi(I)[X,Y] \]

\noindent and so, $\psi(I)=I$ from which the conclusion immediately follows. $\blacksquare$

\end{prf}

Using the scalar matrix $iI$ in place of $I$ in the proof of Proposition (\ref{P5}) we get the following corollary.

\begin{cor}\label{C6}
If the Lie bracket $[,]_{\psi}$, with $\psi(X)=AX$ for all $X\in M_{n}(\mathbb{C})$, coincides with the canonical Lie bracket on $\mathfrak{u}_{S}$ then $A=I$.
\end{cor}

\section{Quantum Channels and Kraus Operators} \label{kraus}

Quantum channels play prominent role in quantum information theory, see \cite{chuang} for more details. They are used to encode operations in the set-up of quantum theory. Quantum channels, in the finite dimensional case, are completely positive maps $M_{n}(\mathbb{C})\stackrel{\Phi}{\longrightarrow}M_{k}(\mathbb{C})$. In some literature, quantum channels are required to be trace-preserving. In this section, we will discuss aspects of quantum channels important for the purpose of this article.

A linear map $A\stackrel{\Phi}{\longrightarrow}B$ between $C^{\ast}$-algebras is said to be $m$-\textit{positive} if the induced map

\[ \Phi_{m}:=I_{m}\otimes\Phi:M_{m}(\mathbb{C})\otimes A\longrightarrow M_{m}(\mathbb{C})\otimes B \]

\noindent sends positive elements to positive elements relative to the natural $C^{\ast}$-algebra structures on the involved tensor products. If $\Phi$ is $m$-positive for all natural numbers $m$ then $\Phi$ is said to be \textit{completely positive}. Note that this notion makes sense for a general map $A\stackrel{\Phi}{\longrightarrow}B$ between $C^{\ast}$-algebras since the minimal and maximal tensor products coincide in the case of tensoring with $M_{n}(\mathbb{C})$. In the case when $A$ and $B$ are finite-dimensional matrix algebras, complete positivity is equivalent to $m$-positivity for some natural number $m$. This is a consequence of Choi's Theorem as stated below. For the proof, see for example \cite{choi} and \cite{mosonyi}.

\begin{thm}{(Choi's Theorem)}\\ \label{cho}
Let $M_{n}(\mathbb{C})\stackrel{\Phi}{\longrightarrow}M_{k}(\mathbb{C})$ be linear. Then, the following are equivalent:
\begin{enumerate}
    \item[(a)] $\Phi$ is completely positive
    \item[(b)] $\Phi$ is $n$-positive
    \item[(c)] There exists $A_{1},\dots,A_{r}\in M_{k,n}(\mathbb{C})$ such that
    \[ \Phi(X)=\sum\limits_{i=1}^{r}A_{i}XA_{i}^{\ast} \]
    
    \noindent for all $X\in M_{n}(\mathbb{C})$. Moreover, the matrices $A_{1},\dots,A_{r}$ satisfy $\sum\limits_{i=1}^{r}A_{i}A_{i}^{\ast}=I$.
\end{enumerate}
\end{thm}

\noindent The representation given in part $(c)$ of the above theorem is called a \textit{Kraus representation} of $\Phi$ and the matrices $A_{1},\dots,A_{r}$ are called the associated \textit{Kraus operators}. The Kraus representation of a linear operator $\Phi$ is far from unique. However, the Kraus representations of a given quantum channel $\Phi$ satisfy certain transitivity relation as stated by the following theorem.

\begin{thm}\label{unitary}
Let $A_{1},\dots,A_{r}$ and $B_{1},\dots,B_{s}$ be two sets of Kraus operators associated to two Kraus representations of a quantum channel $\Phi$. Then, there is a unitary $(U_{ij})\in M_{max\left\{r,s\right\}}(\mathbb{C})$ such that $B_{j}=\sum\limits_{i=1}^{r}U_{ji}A_{i}$.
\end{thm}

\noindent For a proof, see pg. 95 of \cite{mosonyi}.

Linear maps that are not necessarily completely positive have similar representations as that of a Kraus representation of a quantum channel. In the general case, a linear map $M_{n}(\mathbb{C})\stackrel{\Phi}{\longrightarrow}M_{k}(\mathbb{C})$ can be represented as

\[ \Phi(X)=\sum\limits_{i=1}^{r} A_{i}XB_{i}^{\ast} \]

\noindent for some matrices $A_{1},\dots,A_{r},B_{1},\dots,B_{r}\in M_{k,n}(\mathbb{C})$.

\section{Proof of Theorems 1 and 2} \label{proof}

\textsc{Proof of Thereom 1:}
Let $\Phi$ be a linear automorphism extending the Lie algebra isomorphism $\phi$ on the whole $M_{n}(\mathbb{C})$. Let $\Psi$ be its inverse. Then, using the Kraus representation for linear maps $M_{n}(\mathbb{C})\stackrel{\Phi,\Psi}{\longrightarrow}M_{n}(\mathbb{C})$, there are $r$ matrices $A_{i},B_{i}$ and $s$ matrices $A^{\prime}_{j},B^{\prime}_{j}$ such that

\[ \Phi(X)=\sum\limits_{i=1}^{r} A_{i}XB_{i}^{\ast} \hspace{.5in} \text{and} \hspace{.5in} \Psi(X)=\sum\limits_{j=1}^{s} A^{\prime}_{i}X(B^{\prime}_{i})^{\ast} \]

\noindent for all $X\in M_{n}(\mathbb{C})$. Without loss of generality, we can assume the matrices $A_{i}$ and $B_{i}$ form linearly independent sets of matrices. We assume the same for the matrices $A^{\prime}_{i}$ and $B^{\prime}_{i}$. Then, $\Phi\circ \Psi=id$ gives

\[ id(X)=\sum\limits_{i,j} A_{i}A^{\prime}_{j}X(B^{\prime}_{j})^{\ast}B^{\ast}_{i}.\]

\noindent Complete positivity implies that the matrices $A_{i}A^{\prime}_{j}$ and $B_{i}B^{\prime}_{j}$ constitutes a set of Kraus operators for $id$. However, $id(X)=X$ is also a Kraus representation for $id$. Since the Kraus operators appearing in different Kraus representations of the same quantum channel are related by a unitary according to Theorem (\ref{unitary}), we must have $r=s=1$ and so, $\Phi(X)=AXB^{\ast}$. Since $\Phi$ is an isomorphism, the matrices $A$ and $B$ are necessarily invertible. The restriction $\phi$ of $\Phi$ on $\mathfrak{u}_{S}$ is given by $\phi(X)=AXB^{\ast}$. Thus, for any $X,Y\in \mathfrak{u}_{S}$ we have

\[ AXYB^{\ast}-AYXB^{\ast}=\phi[X,Y]=[\phi(X),\phi(Y)]=AXB^{\ast}AYB^{\ast}-AYB^{\ast}AXB^{\ast} \]

\noindent from which we immediately see that

\[ [X,Y]=XY-YX=XB^{\ast}AY-YB^{\ast}AX=[X,Y]_{\psi} \]

\noindent where $\psi(X)=B^{\ast}AX$ for all $X\in\mathfrak{u}_{S}$. Hence, by Corollary (\ref{C6}) we have $B^{\ast}A=I$ and so, $A^{-1}=B^{\ast}$. Taking $V=A$ proves the theorem. $\blacksquare$

\vspace{.25in}

\noindent \textsc{Proof of Theorem 2:}
If $X\in \mathfrak{u}_{S}$ then

\[ T\phi(X)^{\ast}=-\phi(X)T \Longleftrightarrow V^{-1}TV^{-\ast}X^{\ast}=-XV^{-1}TV^{-\ast}. \]

\noindent Thus, $\mathfrak{u}_{S}= \mathfrak{u}_{V^{-1}TV^{-\ast}}$. This implies that the Lie groups $U_{S}$ and $U_{V^{-1}TV^{-\ast}}$ associated to $\mathfrak{u}_{S}$ and $\mathfrak{u}_{V^{-1}TV^{-\ast}}$, respectively, are the same. That is,

\[ SX^{\ast}=X^{-1}S \Longleftrightarrow V^{-1}TV^{-\ast}X^{\ast}=X^{-1}V^{-1}TV^{-\ast} \]

\noindent or equivalently,

\[ XSX^{\ast}=S \Longleftrightarrow (VXV^{-1})T(VXV^{-1})^{\ast}=T \]

\noindent for all $X\in U_{S}$. Thus,

\[ stab(S)=U_{S}=V^{-1}\cdot U_{T}\cdot V=V^{-1}\cdot stab(T)\cdot V \]

\noindent showing that $stab(S)$ and $stab(T)$ are conjugate subgroups of $GL_{n}(\mathbb{C})$. $\blacksquare$

\section{Unanswered Questions} \label{problems}

Relative to the usual commutator bracket $[,]$, the subspaces $\mathfrak{u}_{S}$ are Lie subalgebras of $M_{n}(\mathbb{C})$ for any $S\in M_{n}(\mathbb{C})$. However, as we have seen in Proposition (\ref{P2}), not all $\mathfrak{u}_{S}$ are Lie subalgebras of $M_{n}(\mathbb{C})$ relative to the twisted bracket $[,]_{\psi}$.

\begin{que}
What are the Lie subalgebras of $M_{n}(\mathbb{C})$ under the Lie bracket $[,]_{\psi}$?
\end{que}

The bracket $[,]_{\psi}$ for a linear map $M_{n}(\mathbb{C})\stackrel{\psi}{\longrightarrow}M_{n}(\mathbb{C})$ is a special case of the class of brackets of the form $[X,Y]_{B}:=B(X,Y)-B(Y,X)$ for some bilinear form $B$ on $M_{n}(\mathbb{C})$.

\begin{que}
When is the bracket $[X,Y]_{B}=B(X,Y)-B(Y,X)$ a Lie bracket on $M_{n}(\mathbb{C})$? And on $\mathfrak{u}_{S}$? In particular, when does $[,]_{B}$ satisfy the Jacobi identity?
\end{que}

In Theorem (\ref{T2}), a necessary condition for $\mathfrak{u}_{S}$ and $\mathfrak{u}_{T}$ to be isomorphic is the (unitary) conjugacy of the stabilizers subgroups of the matrices $S$ and $T$. This is much weaker than $S$ and $T$ being unitarily similar.

\begin{que}
Using the orbit-stabilizer theorem for the $\ast$-conjugation action of $GL_{n}(\mathbb{C})$ on $M_{n}(\mathbb{C})$, what can we say about the matrices $S$ and $T$ if $stab(S)=stab(T)$? Or if they are only conjugate subgroups?
\end{que}

The main goal of this article is to determine when the Lie algebras $\mathfrak{u}_{S}$ and $\mathfrak{u}_{T}$ are abstractly isomorphic. Another natural inquiry is to understand the lattice structure of the Lie algebras $\mathfrak{u}_{S}$ in terms of inclusions. In line with this, we have the following question.

\begin{que}
If $\mathfrak{u}_{S}\leqslant\mathfrak{u}_{T}$, is it the case that $stab(T)$ is conjugate to a (possibly trivial) subgroup of $stab(S)$?
\end{que}

More importantly, the author is very much interested with the following question.

\begin{que}
Let $\psi(X)=AX$ for some $A\in M_{n}(\mathbb{C})$. What is the Lie group structure on $GL_{n}(\mathbb{C})$ so that its Lie algebra is $M_{n}(\mathbb{C})$ with bracket $[,]_{\psi}$? By Ado's Theorem, every finite dimensional real Lie algebra $\mathfrak{g}$ is the Lie algebra of a Lie subgroup $G$ of $GL_{R}(\mathbb{C})$ for some $R$. Since the bracket $[,]_{\psi}$ is a 'twist' of the usual commutator on the same vector space, it is reasonable to expect that its Lie group has the same underlying manifold as that of $GL_{n}(\mathbb{C})$ but with a 'twisted' multiplication. See \cite{procesi}.
\end{que}

\subsection*{Acknowledgement}
I would like to thank Diane Pelejo for the valuable communications, bringing to my attention the canonical forms of quantum channels.

\hspace{1in}

\noindent\textsc{Clarisson Rizzie P. Canlubo}\\
University of the Philippines$-$Diliman\\
Quezon City, Philippines 1101\\
crpcanlubo@math.upd.edu.ph


\begin{thebibliography}{1}

  \bibitem{caalim} Caalim, J., Canlubo, C. R. P., and Tanaka, Y.  \textit{On the Lie algebra associated to $S$-unitary matrices.} Linear Algebra and its Applications, 553, pp 167-181, 2018.
  
  \bibitem{choi} Choi, M. -D. \textit{Completely Positive Linear Maps on Complex Matrices}, Linear Algebra and its Applications, 10, pp. 285–290, 1975.
  
  
  \bibitem{chuang} Chuang, C., Nielsen, M.  \textit{Quantum Computation and Quantum Information.} Cambridge University Press, 2010.
  
  \bibitem{gheondea} Gheondea, A.  \textit{The three equivalent forms of completely positive maps on matrices.} Annals of the University of Bucharest (mathematical series), 1 (LIX), 79-98, 2010.
  
  \bibitem{hall} Hall, B.  \textit{Lie Groups, Lie Algebras and Representations. An Elementary Introduction.} Graduate Texts in Mathematics, Springer, 2015.
  
  \bibitem{mosonyi} Mosonyi, M.  \textit{Quantum Information Theory.} (online notes) \url{http://math.bme.hu/~mosonyi/QI/QI_notes171229.pdf}.
  
  \bibitem{nayak} Nayak, A. and Sen, P.  \textit{Invetible Quantum Channels and Perfect Encryption of Quantum States.} \url{https://arxiv.org/abs/quant-ph/0605041}.
  
  \bibitem{procesi} Procesi, C.  \textit{Lie Groups: An Approach through Invariants and Representations.} Universitext, Springer, 2007.
  
  
  
  
\end{thebibliography}
\end{document}